\magnification=1200

\def\Q{{\bf {Q}}}

\def\Z{{\bf Z}}     
\def\R{{\bf R}}

\def\C{{\bf C}}

%%%%%%%%%%%%%%%%%%%%%%%%%%%%%%%%%%%%%%%%%%%%%%%%%%%%
% typoref.tex. V : January 18, 2000. 
% Author : Anthony PHAN
% Warning : syntaxe +- LaTeX 
% Sources :
% T. Lachand--Robert, ``La Ma\^\i trise de \TeX'',
% R\'ef\'erences crois\'ees;
% latex.ltx's sources;
% and of course the \TeX book.
%%%%%%%%%%%%%%%%%%%%%%%%%%%%%%%%%%%%%%%%%%%%%%%%%%%%%
%
\catcode`@=11
%
% style (look at the behavior of \item dans \bibitem too,
% and at one ,\  in \re@dreferenceslist)
% Feel free to change: 	\bibn@me (title like ``R\'ef\'erences'')
%			\bibliographym@rk (general style)
%
\def\bibn@me{R\'ef\'erences}
\def\bibliographym@rk{\centerline{{\sc\bibn@me}}
	\sectionmark\section{\ignorespaces}{\unskip\bibn@me}
	\bigbreak\bgroup
	\ifx\ninepoint\undefined\relax\else\ninepoint\fi}
%
% Beware of the \bgroup: it will be closed by \endthebibliography
%
% \refsp@ce is the spacing command that appens between multiple
% references.
%
\let\refsp@ce=\ 
\let\bibleftm@rk=[
\let\bibrightm@rk=]
%
% if you want more space between brackets...
%\let\refsp@ce=\thinspace
%\def\bibleftm@rk{[\thinspace}
%\def\bibrightm@rk{\thinspace]}
%
% frenchy stuff
%
\def\numero{n\raise.82ex\hbox{$\fam0\scriptscriptstyle o$}~\ignorespaces}
%
% new variables
%
\newcount\equationc@unt
\newcount\bibc@unt
\newif\ifref@changes\ref@changesfalse
\newif\ifpageref@changes\ref@changesfalse
\newif\ifbib@changes\bib@changesfalse
\newif\ifref@undefined\ref@undefinedfalse
\newif\ifpageref@undefined\ref@undefinedfalse
\newif\ifbib@undefined\bib@undefinedfalse
\newwrite\@auxout
%
% mark an equation
%
\def\eqnum{\global\advance\equationc@unt by 1%
\edef\lastref{\number\equationc@unt}%
\eqno{(\lastref)}}
%
% One can reference anything, just copy the former macro
% and use it so: \machin \label{truc}
% In machin you would have defined \lastref by some number
% or any text.
%
% References macros
%
% The next macros are the core of \ref and \cite commands.
% Its first argument may be ref, pageref or bib.
%
% It is too tricky to be explained.
% (It is a bit recursive.)
% It allows using \cite or \ref or ...
% with arbitrary many arguments,
% for instance:
% \cite{knuth1,knuth2,ma pomme}
%
% First argument is always ref, pageref or bib.
%
\def\re@dreferences#1#2{{%
	\re@dreferenceslist{#1}#2,\undefined\@@}}
\def\re@dreferenceslist#1#2,#3\@@{\def\next{#2}%
	\expandafter\ifx\csname#1@@\meaning\next\endcsname\relax
	??\immediate\write16
	{Warning, #1-reference "\next" on page \the\pageno\space
	is undefined.}%
	\global\csname#1@undefinedtrue\endcsname
	\else\csname#1@@\meaning\next\endcsname\fi
	\ifx#3\undefined\relax
	\else,\refsp@ce\re@dreferenceslist{#1}#3\@@\fi}
%
% notice that the former ``,\refsp@ce'' will separate
% multiple arguments. But beware of spaces
% while defining a reference or calling for it!
%
% tricky thing: \newlabel has two arguments
% {labelname}{{\lastref}{\pageref}}
% The second argument is read as two arguments
% by \newl@bel. This was necessary to get
% a jobname.aux containing the same syntax
% LaTeX would produce and use.
%
\def\newlabel#1#2{{\def\next{#1}\newl@bel#2}}
\def\newl@bel#1#2{%
	\expandafter\xdef\csname ref@@\meaning\next\endcsname{#1}%
	\expandafter\xdef\csname pageref@@\meaning\next\endcsname{#2}}
\def\label#1{{%
	\toks0={#1}\message{ref(\lastref) \the\toks0,}%
	\ignorespaces\immediate\write\@auxout%
	{\noexpand\newlabel{\the\toks0}{{\lastref}{\the\pageno}}}%
	\def\next{#1}%
	\expandafter\ifx\csname ref@@\meaning\next\endcsname\lastref%
	\else\global\ref@changestrue\fi%
	\newlabel{#1}{{\lastref}{\the\pageno}}}}
\def\ref#1{\re@dreferences{ref}{#1}}
\def\pageref#1{\re@dreferences{pageref}{#1}}
%
% bibliography macros
%
\def\bibcite#1#2{{\def\next{#1}%
	\expandafter\xdef\csname bib@@\meaning\next\endcsname{#2}}}
\def\cite#1{\bibleftm@rk\re@dreferences{bib}{#1}\bibrightm@rk}
%
% The argument of \beginthebibliography
% is any sequence of numerals which will represent
% the maximum \item's length. If you have less than 9
% \bibitem's, this argument may be {any numeral}.
% if you have between 100 and 999 \bibitem's
% this argument may be {any three numerals},
% and so on.
%
\def\beginthebibliography#1{\bibliographym@rk
	\setbox0\hbox{\bibleftm@rk#1\bibrightm@rk\enspace}
	\parindent=\wd0
	\global\bibc@unt=0
	\def\bibitem##1{\global\advance\bibc@unt by 1
		\edef\lastref{\number\bibc@unt}
		{\toks0={##1}
		\message{bib[\lastref] \the\toks0,}%
		\immediate\write\@auxout
		{\noexpand\bibcite{\the\toks0}{\lastref}}}
		\def\next{##1}%
		\expandafter\ifx
		\csname bib@@\meaning\next\endcsname\lastref
		\else\global\bib@changestrue\fi%
		\bibcite{##1}{\lastref}
		\medbreak
		\item{\hfill\bibleftm@rk\lastref\bibrightm@rk}%
		}
	}
\def\endthebibliography{\egroup\par}
%
% THE NEXT MACRO MUST BE INCLUDED
% IN THE \BYE COMMAND. FOR INSTANCE:
%
% \catcode`@=11
% \outer\def\bye{\@closeaux
% 	\par\vfill\supereject\end}
% \catcode`@=12
%
\def\@closeaux{\closeout\@auxout
	\ifref@changes\immediate\write16
	{Warning, changes in references.}\fi
	\ifpageref@changes\immediate\write16
	{Warning, changes in page references.}\fi
	\ifbib@changes\immediate\write16
	{Warning, changes in bibliography.}\fi
	\ifref@undefined\immediate\write16
	{Warning, references undefined.}\fi
	\ifpageref@undefined\immediate\write16
	{Warning, page references undefined.}\fi
	\ifbib@undefined\immediate\write16
	{Warning, citations undefined.}\fi}
%
% initialization of jobname.aux
%
\immediate\openin\@auxout=\jobname.aux
\ifeof\@auxout \immediate\write16
  {Creating file \jobname.aux}
\immediate\closein\@auxout
\immediate\openout\@auxout=\jobname.aux
\immediate\write\@auxout {\relax}%
\immediate\closeout\@auxout
\else\immediate\closein\@auxout\fi
%
% Let's read this file and open it out
%
\input\jobname.aux
\immediate\openout\@auxout=\jobname.aux
% this file will be closed by \bye.
%
% That's all, folks!
%
\catcode`@=12
%\endinput

%
\catcode`@=11
\def\bibliographym@rk{\bgroup}
%
% \bye est modifie pour la biblio et la table des matieres
%
\outer\def\bye{ 	\par\vfill\supereject\end}

\def\house#1{\setbox1=\hbox{$\,#1\,$}%
\dimen1=\ht1 \advance\dimen1 by 2pt \dimen2=\dp1 \advance\dimen2 by 2pt
\setbox1=\hbox{\vrule height\dimen1 depth\dimen2\box1\vrule}%
\setbox1=\vbox{\hrule\box1}%
\advance\dimen1 by .4pt \ht1=\dimen1
\advance\dimen2 by .4pt \dp1=\dimen2 \box1\relax}

  \def\eps{{\varepsilon}}

\def\build#1_#2^#3{\mathrel{\mathop{\kern 0pt#1}\limits_{#2}^{#3}}}

\def\date {le\ {\the\day}\ \ifcase\month\or janvier
\or fevrier\or mars\or avril\or mai\or juin\or juillet\or
ao\^ut\or septembre\or octobre\or novembre
\or d\'ecembre\fi\ {\oldstyle\the\year}}

\font\fivegoth=eufm5 \font\sevengoth=eufm7 \font\tengoth=eufm10

\newfam\gothfam \scriptscriptfont\gothfam=\fivegoth
\textfont\gothfam=\tengoth \scriptfont\gothfam=\sevengoth

\def\smallsquare{\vbox{\hrule\hbox{\vrule height 1 ex\kern 1 ex\vrule}\hrule}}
\def\cqfd{\hfill \smallsquare\vskip 3mm}

\def\og{\leavevmode\raise.3ex\hbox{$\scriptscriptstyle 
\langle\!\langle\,$}}
\def \fg {\leavevmode\raise.3ex\hbox{$\scriptscriptstyle 
\!\rangle\!\rangle\,\,$}}

\def\rme{{\rm e}}

\def\for{\ \ \hbox{for}\ }

\def\rmh{{\rm h}}

%%%%%%%%%%%%%%%%%%%%%%%%%%%%%%%%%%%%%%%%%%%%
\centerline{}

\vskip 8mm

\centerline{\bf $S$-parts of terms of integer linear recurrence sequences}

\vskip 8mm

\vskip 6mm

\centerline{Yann Bugeaud and Jan-Hendrik Evertse\footnote{}{\rm 
}}

{{\narrower\narrower
\vskip 12mm

\proclaim Abstract. {
Let $S = \{q_1, \ldots , q_s\}$ be a finite, non-empty set of distinct prime numbers.
For a non-zero integer $m$, write $m = q_1^{r_1} \ldots q_s^{r_s} M$, where 
$r_1, \ldots , r_s$ are non-negative integers and $M$ is an integer 
relatively prime to $q_1 \ldots q_s$. We define the $S$-part $[m]_S$ of $m$ by
$[m]_S := q_1^{r_1} \ldots q_s^{r_s}$. 
Let $(u_n)_{n \ge 0}$ be a linear recurrence sequence of integers.
Under certain necessary conditions, 
we establish that for every $\eps > 0$, there exists an integer $n_0$
such that $[u_n]_S\leq |u_n|^{\eps}$ holds for $n > n_0$. 
Our proof is ineffective in the sense that it does not give an explicit value
for $n_0$. Under various assumptions on $(u_n)_{n \ge 0}$, we also give 
effective, but weaker, upper bounds for $[u_n]_S$ of the form $|u_n|^{1 -c}$,
where $c$ is positive and depends only on $(u_n)_{n \ge 0}$ and $S$.
}

}}

\vskip 5mm

\centerline{\bf 1. Introduction and results}

\vskip 5mm

Let $k$ be a positive integer, and let $a_1, \ldots , a_k$ and $u_0, \ldots , u_{k-1}$ be integers
such that $a_k$ is non-zero and $u_0, \ldots , u_{k-1}$ are not all zero. Put
$$
u_n = a_1 u_{n-1} + \ldots + a_{k} u_{n-k}, \quad \hbox{for $n \ge k$}.     \eqno (1.1)
$$
The sequence $(u_n)_{n \ge 0}$ is a linear recurrence sequence of integers of order $k$. Its 
characteristic polynomial 
$$
G(X) := X^k - a_1 X^{k-1} - \ldots - a_k
$$
can be written as 
$$
G(X) = \prod_{i=1}^t \, (X - \alpha_i)^{\ell_i},
$$
where $\alpha_1, \ldots , \alpha_t$ are distinct algebraic numbers and $\ell_1, \ldots , 
\ell_t$ are positive integers. It is then well-known (see e.g. Chapter C in \cite{ShTi86}) 
that there exist 
polynomials $f_1 (X), \ldots , f_t(X)$ of degrees less than $\ell_1, \ldots , \ell_t$,
respectively, and with coefficients in the 
algebraic number field $K := \Q(\alpha_1, \ldots , \alpha_t)$, such that
$$
u_n = f_1(n) \alpha_1^n + \ldots + f_t(n) \alpha_t^n, \quad \hbox{for $n \ge 0$}.   \eqno (1.2)
$$
The recurrence sequence $(u_n)_{n \ge 0}$ is said to be {\it degenerate} if there are 
integers $i, j$ with $1 \le i < j \le t$ such that $\alpha_i / \alpha_j$ is a root of unity.

We keep the above notation throughout the present paper.
The case $t=1$, that is, of sequences $(f(n)a^n)_{n\ge 0}$ where $f(X)$
is an integer polynomial and $a$ a non-zero integer, can be treated using
the work of \cite{GrVi13,BuEvGy17}. 
Thus, in all what follows, 
we assume that $t \ge 2$, the polynomials $f_1 (X), \ldots , f_t(X)$
are non-zero, and $(u_n)_{n \ge 0}$ is non-degenerate. 

By means of a $p$-adic generalization of the Thue--Siegel theorem, 
Mahler \cite{Mah34} proved that every non-degenerate 
binary recurrence sequence $(u_n)_{n \ge 0}$
tends in absolute value to infinity as $n$ tends to infinity. 
This was extended to every non-degenerate 
recurrence sequence by van der Poorten and Schlickewei \cite{vdPSc82} 
and, independently, Evertse \cite{Ev84}. They proved the following stronger result. 
For an integer $m$, let $P[m]$ denote its greatest prime factor, with the 
convention that $P[0]= P[\pm 1] = 1$. By means of a $p$-adic version of the 
Schmidt Subspace Theorem, they established 
that $P[u_n]$ 
%, if the polynomials $f_1(X), \ldots , f_t(X)$ are non-zero, then 
tends to infinity as $n$ tends to infinity.

This result is ineffective, but an effective version of it was proved by
Stewart \cite{Ste82,Ste08c}, under an additional assumption. 
Without loss of generality, assume that
$$
|\alpha_1| \ge |\alpha_2| \ge \ldots \ge |\alpha_t| > 0.
$$
Then, if $|\alpha_1| > |\alpha_2|$ (this assumption is 
often called the {\it dominant root assumption}) 
and $u_n \not= f_1(n) \alpha_1^n$, there are positive numbers $c_1$ and $c_2$, 
which are effectively computable in terms $(u_n)_{n \ge 0}$ (which exactly means,
in terms of $a_1, \ldots , a_k$ and $u_0, \ldots , u_{k-1}$), such that
$$
P[u_n] > c_1 \log n \, {\log \log n \over \log \log \log n}, \quad
\hbox{for $n > c_2$};    \eqno (1.3) 
$$
see also \cite{Ste13a,Ste13b} for stronger results for binary recurrence sequences.

In the present note, we investigate the following related problem. 
Let $S = \{q_1, \ldots , q_s\}$ be a finite, non-empty set of distinct prime numbers.
For a non-zero integer $m$, write $m = q_1^{r_1} \ldots q_s^{r_s} M$, where 
$r_1, \ldots , r_s$ are non-negative integers and $M$ is an integer 
relatively prime to $q_1 \ldots q_s$. We define the $S$-{\it part} $[m]_S$ of $m$ by
%and the $S$-part of $n$ by 
$$
[m]_S := q_1^{r_1} \ldots q_s^{r_s}.
% \quad \hbox{and} \quad (n)_S := |M|.
$$
We ask for a non-trivial upper bound for the $S$-part of the $n$-th term of a 
non-degenerate recurrence sequence of integers.  

A first result on this question was obtained by Mahler \cite{Mah66} in 1966 for a special family
of binary recurrence sequences. 
We keep the above notation and assume that $k=2$, $a_2 \le - 2$, $-4 a_2 > a_1^2$, and
that $a_1$ and $a_2$ are coprime. 
%Let $S$ be a finite, non-empty set of prime numbers and $\eps$ a positive real number. 
By means of a $p$-adic extension of the Roth theorem
established by Ridout \cite{Rid58}, Mahler showed that, if $n$ is large enough, then 
we have $[u_n]_S < |u_n|^{\eps}$. He observed that his result implies that $P[u_n]$ tends to 
infinity as $n$ tends to infinity, a statement which was new at that time.

Before stating our first theorem, which extends Mahler's result to every 
non-degenerate recurrence sequence of integers, we need to introduce
some additional notation. 
Choose embeddings of $K$ in $\C$ and of $K$ in $\Q_p$, for every prime $p$.
These embeddings define extensions to $K$ of the ordinary absolute value $|\cdot |$
and of the $p$-adic absolute value $|\cdot |_p$ for every prime $p$, 
normalized such that $|p|_p = p^{-1}$.
Define the quantity
$$
\delta :=-{\sum_{p\in S}\log\max \{|\alpha_1|_p, \ldots ,  |\alpha_t|_p\} \over
\log\max \{|\alpha_1|, \cdots , |\alpha_t|\} }.    \eqno (1.4) 
$$
By our assumptions on the sequence $(u_n)_{n\geq 0}$, $\alpha_1, \ldots , \alpha_t$
are algebraic integers which are not all roots of unity and whose product
is a non-zero rational integer. Therefore, $\max_i |\alpha_i|>1$,
and $\max_i |\alpha_i|_p\leq 1$ for $p\in S$. 
Hence, $\delta$ is well-defined and $\delta \geq 0$.
Furthermore, we observe that $\delta < 1$ since $t \ge 2$. Indeed, letting 
$A:=\max_i |\alpha_i|$, $A_p:=\max_i|\alpha_i|_p$ for $p\in S$ 
and $a:=\alpha_1\cdots\alpha_t$, we have
$$
(1-\delta )\log A=\log A+\sum_{p\in S} \log A_p
\ge t^{-1}\big(\log |a|+\sum_{p\in S}\log |a|_p\big)\geq 0,
$$
where both inequality signs are equality signs if and only if $|\alpha_1|_p=\cdots =|\alpha_t|_p$
for $p\in S\cup\{\infty\}$ and $|\alpha_1|_p=\cdots =|\alpha_t|_p=1$ 
for all prime numbers $p$ outside $S$,
that is, if all quotients $\alpha_i/\alpha_j$ are roots of unity, which is against
our assumption. So the left-hand side of the above inequality is $>0$, thus $\delta < 1$. 

Our first result is an easy consequence of work of Evertse,
see e.g., \cite{Ev84}, Theorem 2, or Proposition 6.2.1 of
\cite{EvGy15}.

\proclaim Theorem 1.1. 
Let $(u_n)_{n \ge 0}$ be a non-degenerate recurrence sequence of integers defined in (1.1). 
Let  $S := \{q_1, \ldots , q_s\}$ be a finite, non-empty set of prime numbers, and $\delta$
be as in (1.4). Further, let $\eps >0$. Then
for every sufficiently large $n$ we have
$$
|u_n|^{\delta -\eps}\leq [u_n]_S\leq |u_n|^{\delta +\eps}.    \eqno (1.5)
$$ 
In particular, if $\gcd (q_1\cdots q_s,a_1, \ldots ,  a_k)=1$, we have for every sufficiently
large $n$,
$$
[u_n]_S\leq |u_n|^{\eps}.
$$

Observe that the assumption $\gcd (q_1\cdots q_s,a_1, \ldots ,  a_k)=1$
implies that $\delta =0$.
Indeed, suppose $\delta >0$. Then, there is a prime number $p$ in $S$ such that
$\max_i |\alpha_i|_p<1$. Since $a_1, \ldots , a_k$ 
are up to sign the elementary symmetric functions
in the $\alpha_i$ taken $\ell_i$ times, we must then have
$|a_i|_p<1$ for $i=1, \ldots ,  k$.
This is clearly impossible.

The proof of Theorem 1.1 depends ultimately on the $p$-adic Schmidt Subspace Theorem,
therefore we cannot compute the set of $n$ for which (1.5) does not hold. 
It seems to be difficult to bound its 
cardinality, even by means of the strongest available versions of the quantitative 
Subspace Theorem.

Our main effective theorem is the following.

\proclaim Theorem 1.2. 
Let $(u_n)_{n \ge 0}$ be a non-degenerate recurrence sequence of integers 
having a dominant root. 
Let $S$ be a finite, non-empty set of prime numbers. 
Then, there exist effectively computable positive numbers $c_1$ and $c_2$, depending 
only on $(u_n)_{n \ge 0}$ and $S$, such that
$$
[u_n]_S \le |u_n|^{1 - c_1}, \quad
%\hbox{for every $n \ge c_2$}.
$$
for every $n \ge c_2$.

Removing the dominant root assumption seems to be very difficult. However, this can be 
done for non-degenerate binary recurrence sequences of integers.

\proclaim Theorem 1.3. 
Let $(u_n)_{n \ge 0}$ be a non-degenerate binary recurrence sequence of integers. 
Assume that $u_n = a \alpha^n + b \beta^n$ for $n \ge 0$, with $a b \alpha \beta \not= 0$. 
Let $S$ be a finite, non-empty set of prime numbers. 
Then, there exist effectively computable positive numbers $c_1$, 
depending only on $(u_n)_{n \ge 0}$, and $c_2$, depending 
only on $(u_n)_{n \ge 0}$ and $S$, such that
$$
[u_n]_S \le |u_n|^{1 - c_1}, \quad
\hbox{for every $n \ge c_2$}.
$$

We stress that, in Theorem 1.3, 
the number $c_1$ is independent of the set $S$ of prime numbers. 
Clearly, when $\alpha$ and $\beta$ are complex conjugates, there is no 
dominant root.

Our next statement is a $p$-adic analogue to Theorem 1.2. 

Let $p$ be a prime number. 
We say that the recurrence sequence $(u_n)_{n \ge 0}$ as in (1.2) 
has a $p$-adic dominant root if there exists $j$ such that $1 \le j \le t$ and 
$|\alpha_j|_p = 1$, while $|\alpha_i|_p < 1$ for $1 \le i \le t$ with $i \not= j$.

\proclaim Theorem 1.4. 
Let $p$ be a prime number.
Let $(u_n)_{n \ge 0}$ be a non-degenerate recurrence sequence of integers 
having a $p$-adic dominant root. 
Let $S$ be a finite, non-empty set of prime numbers. 
Then, there exist positive numbers $c_1$ and $c_2$, depending 
only on $(u_n)_{n \ge 0}$ and $S$, such that
$$
[u_n]_S \le |u_n|^{1 - c_1}, \quad
$$
for every $n \ge c_2$. Furthermore, for every positive real number $\eps$, there 
exists an effectively computable integer $c_3$, depending only on 
$(u_n)_{n \ge 0}$ and $\eps$, such that 
$$
P[u_n] > (1 - \eps) \, \log n \, {\log \log n \over \log \log \log n}, \quad
\hbox{for $n > c_3$}.    
$$

The last statement of Theorem 1.4 seems to be new.

The proof of Theorem 1.2 allows us to establish the following statement.

\proclaim Theorem 1.5. 
Let $\theta > 1$ be a real algebraic number such that all of its Galois conjugates 
are less than $\theta$ in modulus. 
Let $\lambda$ be a non-zero real algebraic number. 
Let $S$ be a finite set of prime numbers. 
If $\theta^\ell\not\in{\bf Z}$ for every integer $\ell \ge 1$,
then there exist effectively computable positive numbers $c_1$ and $c_2$, depending 
only on $\lambda$, $\theta$, and $S$, such that
$$
[ \lfloor \lambda \theta^n \rfloor ]_S \le |\lambda \theta^n|^{1 - c_1}, \quad
\hbox{for every $n \ge c_2$}.
$$

Theorem 1.5 applies to the sequence of integer parts of $(3/2)^n$. 
Lower bounds for the greatest prime factor of $\lfloor \theta^n \rfloor$, where 
$\theta > 1$ is an algebraic number such that $\theta^\ell$ is not an integer for 
every integer $\ell \ge 1$, have been obtained by Luca and Mignotte \cite{LuMi09}. 
They are similar to (1.3). 

It would be interesting to see under which assumption on 
the algebraic numbers $\lambda$ and $\theta$ we get 
$$
[ \lfloor \lambda \theta^n \rfloor ]_S \le  |\lambda \theta^n|^\eps,
$$
for every $\eps > 0$ and every sufficiently large integer $n$.

\vskip 5mm

\goodbreak

\centerline{\bf 2. Proof of Theorem 1.1}

\vskip 5mm

Recall that we have set $K:=\Q (\alpha_1, \ldots , \alpha_t)$. 
Denote by $M_K$ the set of places
of $K$. For $v$ in $M_K$, we choose a normalized absolute value $|\cdot |_v$
such that if $v$ is an infinite place, then 
$$
|x|_v=|x|^{[K_v:\R ]/[K:\Q ]}, \for x\in\Q,
$$
while if $v$ is finite and lies above the prime $p$, then
$$
|x|_v=|x|_p^{[K_v:\Q_p ]/[K:\Q ]}, \for x\in\Q. 
$$
These absolute values satisfy
the product formula 
$$
\prod_{v\in M_K} |x|_v=1, \quad \hbox{for every non-zero $x\in K$}.
$$
Moreover, if $x\in\Q$, then $\prod_{v|\infty} |x|_v=|x|$ 
and $\prod_{v|p} |x|_v=|x|_p$, where the products are taken
over all infinite places of $K$, respectively all places of $K$
lying above the prime number~$p$.

Let $T$ be a finite set of places of $K$, containing all infinite places.
Define the ring of $T$-integers and the group of $T$-units of $K$ by
$$
O_T:=\{ x\in K:\, |x|_v\leq 1\for v\in M_K\setminus T\},
$$
$$
O_T^*:=\{ x\in K:\, |x|_v= 1\for v\in M_K\setminus T\},
$$
respectively. Further define
$$
H_T(x_1, \ldots ,  x_n):=\prod_{v\in T}\max \{|x_1|_v, \ldots ,  |x_n|_v\}, 
\for x_1, \ldots ,  x_n\in O_T.
$$
Our main tool is the following result, which is Proposition 6.2.1 of \cite{EvGy15}
and which is essentially the same as Theorem 2 of \cite{Ev84}.

\proclaim Proposition 2.1. 
Let $U$ be a subset of $T$, $t\geq 2$ and $\eps >0$. Then for all
$x_1, \ldots ,  x_t\in O_T$ such that every non-empty subsum of $x_1+\cdots +x_t$ 
is non-zero, we have
$$
\prod_{i=1}^t\prod_{v\in T} |x_i|_v\cdot \prod_{v\in U} |x_1+\cdots +x_t|_v
\gg \prod_{v\in U} \max \{|x_1|_v, \ldots ,  |x_t|_v\} \cdot H_T(x_1, \ldots ,  x_t)^{-\eps},
$$
where the implied constant depends on $K,T,t$ and $\eps$.

The proof of this result depends on the $p$-adic Schmidt Subspace Theorem,
therefore the implied constant is ineffective.

\vskip0.3cm\noindent
{\it Proof of Theorem 1.1.}
We introduce the following notation. First, by $O(1)$ we denote constants
depending on $(u_n)_{n\geq 0}$ and $S$. Second, we choose a real number 
$\eps'>0$ which will later be taken sufficiently small in terms
of $\eps$; then constants implied by the Vinogradov symbols $\ll$, $\gg$
will depend on $(u_n)_{n\geq 0}$, $S$ and $\eps'$.
Lastly, we put
$$
A:=\max \{|\alpha_1|, \ldots ,  |\alpha_t|\},\ \ 
A_p:=\max \{|\alpha_1|_p, \ldots ,  |\alpha_t|_p\}\, \for p\in S,
$$
and
$$
A_v:=\max \{|\alpha_1|_v, \ldots ,  |\alpha_t|_v\}\, \for v\in M_K .
$$
Then by our choice of the absolute values on $K$, we have
$$
\prod_{v|\infty} A_v=A,\ \ \ \prod_{v|p} A_v=A_p\, \for p\in S.
$$ 

We choose a finite set of places $T$ of $K$, containing all infinite places
and all places lying above the primes in $S$, 
such that $\alpha_1, \ldots ,  \alpha_t\in O_T^*$ and the coefficients
of $f_1(X), \ldots ,  f_t(X)$ are in $O_T$.
Let for the moment $U$ be any subset of $T$.
Each subsum of $u_n=\sum_{i=1}^t f_i(n)\alpha_i^n$ is a non-degenerate
linear recurrence sequence of algebraic numbers in $K$. 
So by the Skolem-Mahler-Lech Theorem,
there are only finitely many non-negative integers $n$ for which at least one 
of the subsums of $u_n$ vanishes. Then by Proposition 2.1 we have
for the remaining positive integers $n$,
$$
\prod_{i=1}^t\prod_{v\in T}|f_i(n)\alpha_i^n|_v\cdot\prod_{v\in U} |u_n|_v
\gg\prod_{v\in U}\max_{1\leq i\leq t} |f_i(n)\alpha_i^n|_v
\cdot \Big(\prod_{v\in T}\max_{1\leq i\leq t} |f_i(n)\alpha_i^n|_v\Big)^{-\eps'/2}.
$$
Since $\prod_{v\in T} |\alpha_i|_v=1$ for $i=1, \ldots ,  t$,
the left-hand side of this inequality is 
$$
\ll (2n)^{O(1)}\prod_{v\in U} |u_n|_v, 
$$
while the right-hand side is 
$$
\gg (2n)^{-O(1)}\Big(\prod_{v\in U} A_v\Big)^n
\cdot\Big(\prod_{v\in T} A_v\Big)^{-n\eps'/2}\gg (2n)^{-O(1)} 
\Big(\prod_{v\in U} A_v\Big)^n\cdot A^{-n\eps'/2},
$$
where we have used $A_v\leq 1$ if $v$ is finite and $\prod_{v|\infty} A_v=A$.
Thus,
$$
\prod_{v\in U}|u_n|_v\gg (2n)^{-O(1)} \Big(\prod_{v\in U} A_v\Big)^n\cdot A^{-n\eps'/2}.
$$
On the other hand, we have a trivial upper bound
$$
\prod_{v\in U}|u_n|_v\leq (2n)^{O(1)}\Big(\prod_{v\in U} A_v\Big)^n.
$$
Since $A>1$, this implies that for every sufficiently large $n$,
$$
\Big(\prod_{v\in U} A_v\Big)^n\cdot A^{-n\eps'}\leq
\prod_{v\in U}|u_n|_v\leq \Big(\prod_{v\in U} A_v\Big)^n\cdot A^{n\eps'}.
$$
We apply this with two choices of $U$.
First let $U$ consist of the infinite places of $K$.
Then, 
$$
\prod_{v\in U} |u_n|_v=|u_n| \quad \hbox{and} \quad  \prod_{v\in U} A_v=A, 
$$
implying that for all sufficiently large $n$,
$$
A^{n(1-\eps')}\leq |u_n|\leq A^{n(1+\eps')}.
$$
Next, let $U$ consist of the places of $K$ lying above the primes in $S$.
Then 
$$
\prod_{v\in U} |u_n|_v=\prod_{p\in S}|u_n|_p=[u_n]_S^{-1}
$$
and 
$$
\prod_{v\in U} A_v=\prod_{p\in S} A_p=A^{-\delta}, 
$$
hence
$$
A^{n(\delta -\eps')}\leq [u_n]_S\leq A^{n(\delta +\eps')}
$$
for every sufficiently large $n$. By taking $\eps'$ sufficiently small
in terms of $\eps$, Theorem 1.1 easily follows.
\cqfd

\vskip 5mm

\centerline{\bf 3. Proofs of Theorems 1.2 and 1.3}

\vskip 5mm

As usual, $\rmh (\alpha)$ denotes the (logarithmic) Weil height of the 
algebraic number $\alpha$. 

Our first auxiliary result is an immediate 
corollary of a theorem of Matveev \cite{Matv00}.

\proclaim Theorem 3.1. 
Let $n \ge 2$ be an integer, 
let $\alpha_1, \ldots, \alpha_n$ be non-zero algebraic numbers and 
let $b_1, \ldots, b_n$ be integers. 
Further, let $D$ be the degree over $\Q$ 
of a number field containing the $\alpha_i$, and
let $A_1, \ldots, A_n$ be real numbers with
$$
\log A_i
\ge
\max
\Bigl\{\rmh (\alpha_i), {|\log \alpha_i| \over D}, {0.16 \over D}
\Bigr\},
\qquad
1\le i \le n.
$$
Set
$$
B := \max\Bigl\{1, \max\Bigl\{ |b_j| \ {\log A_j \over \log A_n} : 1 \le j \le n \Bigr\} \Bigr\}.
$$
Then, we have
$$
\log |\alpha_1^{b_1} \ldots \alpha_n^{b_n} - 1|  
> - 4 \times 30^{n+4} \, (n+1)^{5.5} \, D^{n+2} \, \log (\rme D) \,  \log (\rme n B) \, 
\log A_1 \ldots \log A_n.      
$$

The key point for our main theorem is the factor $\log A_n$ in the denominator in the
definition of $B$. 

Our second auxiliary result is extracted from \cite{Yu07}. 

\proclaim Theorem 3.2. 
Let $p$ be a prime number, $K$ an algebraic number field of degree $D$
and $|\cdot |_p$ an absolute value on $K$ with $|p|_p=p^{-1}$.
Further,
let $\alpha_1, \ldots, \alpha_n$ be elements of $K$,
and let $A_1, \ldots , A_n$ be real numbers with
$$
\log A_i \ge
\max\Bigl\{ \rmh(\alpha_i), {1 \over 16 \rme^2 D^2} \Bigr\},
\qquad 1\le i \le n.
$$
Let $b_1, \ldots , b_n$ denote nonzero rational integers and 
let $B$ and $B_n$ be real numbers such that
$$
B \ge \max\{|b_1|, \ldots , |b_n|, 3\} \quad 
\hbox{and} \quad
B \ge B_n \ge |b_n|.
$$
Assume that
$$
%v_p (b_n) \le v_p (b_j), \quad j = 1, \ldots , n.
|b_n|_p \ge |b_j|_p,  \quad j = 1, \ldots , n.
$$
Let $\delta$ be a real number with $0< \delta \le 1/2$. 
With the above notation, we have
$$
\eqalign{
%v_p (\alpha_1^{b_1} \ldots \alpha_n^{b_n} - 1)  < &  \, 
\log |\alpha_1^{b_1} \ldots \alpha_n^{b_n} - 1|_p  > &  \, -
(16 \rme D)^{2(n+1)} n^{3/2}  \, (\log (2nD))^2
 \, D^n {p^D \over \log p} \times  \cr
&  \, \, \, \, \, \, \times 
\max\Bigl\{  (\log A_1) \cdots (\log A_n) (\log T), {\delta B \over B_n c_0(n, D)} \Bigr\}, \cr}
$$
where
$$
T = B_n \delta^{-1} c_1(n, D) p^{(n+1)D}  (\log A_1) \cdots (\log A_{n-1})
$$
and
$$
c_0(n, D) = (2D)^{2n+1} \log (2D) \log^3 (3D), \quad
c_1(n, D) = 2 \rme^{(n+1)(6n+5)} D^{3n} \log (2D).
$$

\bigskip

\noindent {\it Proof of Theorem 1.2.}  
We establish a slightly more general result.
 
We consider a sequence of integers $(v_n)_{n \ge 0}$ 
with the property that there are $\theta$ in $(0, 1)$ and $C > 0$ such that
$$
| v_n - f(n) \alpha^n | \le C \, |\alpha|^{\theta n},  \quad n \ge 0,     \eqno (3.1) 
$$
where $f (X)$ is a non-zero polynomial whose coefficients are algebraic numbers
and $\alpha$ is an algebraic number with $|\alpha|>1$
Clearly, a recurrence sequence having a dominant root $\alpha$ has the above property.
We prove a similar statement as Theorem 1.2 for the sequence $(v_n)_{n\ge 0}$.  

The constants $c_1, c_2, \ldots$ below are positive, effectively computable and depend 
at most on $(v_n)_{n \ge 0}$. 
The constants $C_1, C_2, \ldots$ below are positive, effectively computable and absolute.

Let $q_1, q_2, \ldots , q_s$ be distinct prime numbers written in increasing order. 
Let $n$ be a positive integer such that $f(n) v_n$ is non-zero and $v_n \not= f(n) \alpha^n$. 
There exist non-negative integers $r_1, \ldots , r_s$ and a non-zero integer $M$ 
coprime with $q_1 \ldots q_s$ such that
$$
v_n = q_1^{r_1} \cdots q_s^{r_s} M.
$$
Observe that there exist positive real numbers $c_1$ and $c_2$ such that
$$
r_j \log q_j \le c_1 n, \quad j = 1, \ldots , s,     \eqno (3.2)
$$
and
$$
\Lambda := |q_1^{r_1} \cdots q_s^{r_s} (M f(n)^{-1}) \alpha^{-n}  - 1| 
\le c_2 |f(n)|^{-1} \, |\alpha|^{(\theta -1) n}. 
$$
Since $\theta < 1$, we get by (3.1) that 
$$
\log |\Lambda | \le - c_3 \,  n. 
$$
Setting
$$
Q := (\log q_1) \cdots (\log q_s) \quad 
\hbox{and} \quad   \log A := \max\{\rmh (M f(n)^{-1}), 2\}, 
$$
Theorem 3.1 and (3.2) imply that
$$
\log |\Lambda | \ge - c_4 C_1^s \, Q \, (\log A) \, \log {n  \over \log A}.
$$
%where $c_1$ is a positive real number.
%Since $Q$ is at least equal to the product of the $s$ first prime numbers, we derive
Comparing both estimates, we obtain that
$$
n  \le c_5 \,  C_2^s \, Q \, (\log Q) \, (\log A).    \eqno (3.3)
$$
Observe that
$$
\log A \le \log |M| + c_6 \log n. 
$$
%for some positive $c_4$ depending only on ${\bf v}$. 
We distinguish two cases.

If $|M| \ge n^{c_6}$, then $A \le M^2$, and, by (3.3), 
$$
n  \le c_7 \,  C_2^s \, Q \, (\log Q) \, (\log |M|).     
$$
We derive that
$$
{|v_n|\over [v_n]_S} = |M| \ge  2^{c_8 n (C_2^s \, Q \,  (\log Q))^{-1}}
\ge |v_n|^{c_9 (C_2^s \, Q \,  (\log Q))^{-1}},
$$
since $|v_m| \le |\alpha_1|^{2m}$ for $m$ sufficiently large. 

If $|M| < n^{c_6}$, then $\log A \le 2 c_6 \log n$ and, by (3.3), 
$$
n  \le c_{10} \,  C_2^s  \, Q \, (\log Q) \, (\log n),
$$
thus, 
$$
n   \le c_{11} \, C_3^s \, Q \,  (\log Q)^2.     \eqno (3.4) 
$$
This completes the proof of Theorem 1.2. \cqfd 

\medskip

\noindent {\bf Remark. } 
Let $\eps$ be a positive real number. 
In the particular case where $M = \pm 1$ and $q_1, \ldots , q_s$ 
are the first $s$ prime numbers $p_1, \ldots , p_s$, we get
from (3.4) and the Prime Number Theorem that
$$
\log n \le c_{12} + s C_4 + (1 + \eps) \, \sum_{k=1}^s \log \log p_k  
\le (1+ 2 \eps) p_s {\log \log p_s \over \log p_s},
$$
if $n$ is sufficiently large in terms of $\eps$. 
This gives
$$
P[v_n]  \ge (1 - 3 \eps) \,  \log n \, {\log \log n \over \log \log \log n}, 
$$
if $n$ is sufficiently large in terms of $\eps$, and we recover Stewart's result (1.3). 
\cqfd

\vskip 4mm

\bigskip

\noindent {\it Proof of Theorem 1.3.}  
The constants $c_1, c_2, \ldots$ below are positive, effectively computable and depend 
at most on $(u_n)_{n\ge 0}$.

Let $q_1, q_2, \ldots , q_s$ be distinct prime numbers written in increasing order. 
Let $n \ge 2$ be an integer. 
There exist non-negative integers $r_1, \ldots , r_s$ and a non-zero integer $M$ 
coprime with $q_1 \ldots q_s$ such that
$$
u_n = q_1^{r_1} \cdots q_s^{r_s} M.
$$
It follows from inequality (20) of \cite{Ste82} that for $i=1, \ldots , s$ we have
$$
r_i \le c_1 \, {q_i^2 \over \log q_i} \, (\log n)^2.
$$
By Lemma 6 of \cite{Ste82}, we get
$$
\log |u_n| > c_2 n.
$$
Consequently, setting
$$
Q := \sum_{i=1}^s \,  q_i^2, 
$$
we obtain
$$
c_2 n < \log |M| + c_1 \, Q \, (\log n)^2.
$$
We distinguish two cases.

If $\log |M| > c_1 \, Q \, (\log n)^2$, then
$
c_2 n <  2 \log |M|,
$
thus
$$
|M| > |u_n|^{c_3}.
$$
If $\log |M|  \le c_1 \, Q \, (\log n)^2$, then
$$
n < c_4 \, Q \, (\log n)^2,
$$
and $n$ is bounded in terms of $a, b, \alpha, \beta$, and $S$. 
This completes the proof of the theorem.  
\cqfd

\bigskip

\noindent {\it Proof of Theorem 1.4.}  
We establish a slightly more general result.
Let $p$ be a prime number, $K$ a number field and $|\cdot |_p$
an absolute value on $K$ with $|p|_p=p^{-1}$.  
We consider a sequence of integers $(v_n)_{n \ge 0}$ 
with the property that there are $\theta$ in $(0, 1)$ and $C > 0$ such that
$$
| v_n - f(n) \alpha^n |_p \le C \, p^{- \theta n},  \quad n \ge 0,   \eqno (3.5) 
$$
where $f (X)$ is a non-zero polynomial with coefficients in $K$
and $\alpha$ an element of $K$ with $|\alpha|_p=1$. 
We also assume that there exist $\beta > 1$ and $n_0$ such that 
$$
v_n\not= f(n)\alpha^n,\ \ 
|v_n| <   \beta^n,  \quad \hbox{for every integer $n > n_0$.}   \eqno (3.6)
$$
Clearly, every non-degenerate recurrence sequence having a $p$-adic
dominant root satisfies both properties.
We prove a statement analogous to Theorem 1.4 for the sequence $(v_n)_{n\ge 0}$. 

Let $q_1, q_2, \ldots , q_s$ be distinct prime numbers written in increasing order. 
Let $n \ge 3$ be an integer such that $f(n) v_n$ is non-zero and $v_n \not= f(n) \alpha^n$. 
There exist non-negative integers $r_1, \ldots , r_s$ and a non-zero integer $M$ 
coprime with $q_1 \ldots q_s$ such that
$$
v_n = q_1^{r_1} \cdots q_s^{r_s} M.
$$

The constants $c_1, c_2, \ldots$ below are positive, effectively computable and depend 
at most on $(v_n)_{n \ge 0}$ and $p$. 
The constants $C_1, C_2, \ldots$ below are positive, effectively computable and absolute.

Consider the quantity
$$
\Lambda := q_1^{r_1} \cdots q_s^{r_s} (M  f(n)^{-1}) \alpha^{-n} - 1
$$
and observe that, by (3.5), we have
$$
\log |\Lambda|_p < - c_1 n. 
$$
Set
$$
Q := (\log q_1) \cdots (\log q_s) \quad 
\hbox{and} \quad   \log A := \max\{\rmh (M f(n)^{-1}), 2\}. 
$$
It follows from Theorem 3.2 applied with 
$$
B:= \max\{ r_1, \ldots , r_s, n\} \quad \hbox{and} \quad \delta = {Q (\log A) \over B}
$$
that
$$
B < 2 Q \log A, \quad \hbox{if $\delta > 1/2$},   \eqno (3.7)
$$
and, otherwise,
$$
\log |\Lambda|_p > - c_1 C_1^s Q (\log A) \max\Bigl\{ \log \Bigl( {B \over \log A} \Bigr), 1 \Bigr\}. 
\eqno (3.8)
$$
If (3.7) holds, then we have
$$
n \le B < 2 Q \log A.    \eqno (3.9) 
$$
If (3.8) holds, then, using
$$
B \le c_3 n,   
$$
which follows from (3.6), we obtain an upper bound similar to (3.3), namely
$$
n \le c_4 C_2^s Q (\log Q) (\log A).    \eqno (3.10)
% \le c_8 C_3^s Q (\log Q) (\log |M| + c_9 \log n),
$$
In view of (3.9) and (3.10), we conclude as in the proof of Theorem 1.2. 
The last statement of the theorem corresponds to the remark following the 
proof of Theorem 1.2. 
\cqfd

\vskip 7mm

\centerline{\bf References}

\beginthebibliography{999}

\medskip 

\bibitem{BuEvGy17}
Y. Bugeaud, J.-H. Evertse, and K. Gy\H ory,
{\it $S$-parts of values of polynomials and of decomposable forms at integral points}.
Preprint. 

\bibitem{Ev84}
J.-H. Evertse,
{\it On sums of $S$-units and linear recurrences}, 
Compos. Math. 53 (1984), 225--244.

\bibitem{EvGy15}
J.-H. Evertse and K. Gy\H ory,
Unit Equations in Diophantine Number Theory.
Cambridge University Press, 2015.

\bibitem{GrVi13}
S. S. Gross and A. F. Vincent,
{\it On the factorization of $f(n)$ for $f(x)$ in $\Z[x]$}, 
Int. J. Number Theory 9 (2013), 1225--1236.

\bibitem{LuMi09}
F. Luca and M. Mignotte,
{\it Arithmetic properties of the integer part of the powers of an algebraic number}, 
Glas. Mat. Ser. III 44 (2009), 285--307.

\bibitem{Mah34}
K. Mahler,
{\it Eine arithmetische Eigenschaft der rekurrierenden Reihen},
Mathematica (Zutphen) 3 (1934), 153--156.

\bibitem{Mah66}
K. Mahler,
{\it A remark on recursive sequences},
J. Math. Sci. Delhi 1 (1966), 12--17.

\bibitem{Matv00} 
E.\ M.\ Matveev,
{\it An explicit lower bound for a homogeneous rational linear form
in logarithms of algebraic numbers.\ II},
Izv.\ Ross.\ Acad.\ Nauk Ser.\ Mat.\  {64}  (2000),  125--180 (in Russian); 
English translation in Izv.\ Math.\  {64} (2000),  1217--1269.

\bibitem{vdPSc82}
A. J. van der Poorten and H. P. Schlickewei, 
{\it The growth condition for recurrence sequences}, 
Macquarie University Math. Rep. 82-0041 (1982). 

\bibitem{Rid58}
D. Ridout,
{\it The p-adic generalization of the Thue-Siegel-Roth theorem}, 
Mathematika 5 (1958), 40--48.

\bibitem{ShTi86}
T. N. Shorey and R. Tijdeman,  
Exponential Diophantine equations.
Cambridge Tracts in Mathematics, 87. Cambridge University Press, Cambridge, 1986.

\bibitem{Ste82}
C. L. Stewart,
{\it On divisors of terms of linear recurrence sequences},
J. reine angew. Math.  333  (1982), 12--31.

\bibitem{Ste08c}
C. L. Stewart,
{\it On the greatest square-free factor of terms of a linear recurrence sequence}. 
In: Diophantine equations, 257--264, 
Tata Inst. Fund. Res. Stud. Math., 20, Tata Inst. Fund. Res., Mumbai, 2008.

\bibitem{Ste13a}
C. L. Stewart,
{\it On divisors of Lucas and Lehmer numbers},
Acta Math. 211 (2013),  291--314. 

\bibitem{Ste13b}
C. L. Stewart,
{\it On prime factors of terms of linear recurrence sequences}. 
In: Number theory and related fields, 341--359, 
Springer Proc. Math. Stat., 43, Springer, New York, 2013.

\bibitem{Yu07}
K. Yu, 
{\it $p$-adic logarithmic forms and group varieties. III}, 
Forum Math. 19 (2007), \ \ 187--280.

\vskip1cm

\vbox{
\hbox{Yann Bugeaud \hfill}
\hbox{Universit\'e de Strasbourg, CNRS  \hfill}
\hbox{IRMA UMR 7501 \hfill}
\hbox{7, rue Ren\'e Descartes \hfill}
\hbox{67000 STRASBOURG (FRANCE)\hfill}
\vskip2pt\noindent
\hbox{{\tt bugeaud@math.u-strasbg.fr} \hfill}
}
\vskip1cm\noindent
\vbox{
\hbox{Jan-Hendrik Evertse \hfill}
\hbox{Universiteit Leiden \hfill}
\hbox{Mathematisch Instituut \hfill}
\hbox{Postbus 9512 \hfill}
\hbox{2300 RA LEIDEN (THE NETHERLANDS)\hfill}
\vskip2pt\noindent
\hbox{{\tt evertse@math.leidenuniv.nl} \hfill}
}

\bye